\documentstyle[12pt]{article} 
\begin{document} 
{ \pagestyle{empty} 
\vspace{10mm} 
\centerline{\Large \bf Metric derived from Lie Groups}
\vspace{10mm}
\centerline{Kazuyasu Shigemoto \footnote{E-mail address:
shigemot@tezukayama-u.ac.jp}}
\centerline {{\it Tezukayama University, Nara 631, Japan }}
\vspace{10mm}

\centerline{\bf Abstract}

We give the expression of the metric derived from Lie groups.
For the metric derived from classical Lie groups such as the 
unitary group, the orthogonal group and the symplectic group, 
we conjecture that the metric becomes the Einstein metric. 
}

\vspace{20mm}

\setcounter{equation}{0}
\section{Introduction} 
In the theory of general relativity, the symmetry of the metric is quite 
important to classify the space-time. 
If we start from the metric itself, it is not easy to  
to find the symmetry of the metric. 
Then we take another approach, where we can find the symmetry of the 
metric in a trivial way. We start from a certain Lie group 
\cite{Jacobson,Hamermesh,Georgi,Gilmore,Cahn} and 
derived the metric from this Lie group. Then the symmetry
of the metric is that of the Lie group. 

The metric derived 
from Lie group has the quite simple structure. The metric 
on $n$-sphere $S^n$ is known to have the Einstein metric
$R_{\mu \nu}=2\Lambda g_{\mu \nu}$ \cite{Besse}. The symmetry of 
$S^n$ is the subset of the symmetry of $SO(n)$, then we expect that
the metric derived from $SO(n)$ becomes the Einstein metric.
In this paper, we derived the metric from $SU(2)$ Lie group in 
two different ways. We then confirm that the metric derived from 
$SU(2)$ becomes the Einstein metric. For classical groups such as 
$SO(n)$, $SU(n)$, $Sp(n)$, the metric derived from the Lie group is 
connected with the metric on the manifold which is the quadratic 
invariant form to define Lie groups. From this observation, 
we conjecture that the metric derived from the classical groups 
becomes the Einstein metric.

\vspace{20mm}

\section{ Metric derived from Lie groups}
In this section, we first give the general method to derive  
the metric from Lie groups. Next we explicitly give the expression
of the metric derived from the $SU(2)$ Lie group in two different ways.

\subsection{General method to derive the metric from Lie groups}
\indent
We denote $U(\theta_1,\theta_2, \cdots)$ as the representation 
of any Lie group, which is parametrized by $\{ \theta_1,
\theta_2, \cdots \}$.  
The global Lie group transformation is given by 
\begin{eqnarray}
    U'=VU   ,
\label{e1}
\end{eqnarray}
where $V$ is the global Lie group element.
In order to construct the metric, we define the following line element
\begin{eqnarray}
    (ds)^2=k {\rm Tr} \left(\overline{U^{-1} dU}  U^{-1} dU 
    \right)   .
\label{e2}
\end{eqnarray}
This line element $(ds)^2$ is the desired line element for 
the Lie groups. The right-hand side of Eq.(\ref{e2}) is 
invariant under the transformation Eq.(\ref{e1}), because we have 
\begin{eqnarray}
   U^{-1}dU \rightarrow U^{-1}V^{-1} V dU =U^{-1} dU .
\label{e3}
\end{eqnarray} 
From the expression of Eq.(\ref{e2}), the above line 
element $(ds)^2$  becomes real number. Then the above line
element satisfy the physically meaningful conditions, that is , 
$i)$ it is the quadratic form of the infinitesimal quantity, 
$ii)$ it is real, $iii)$ it is invariant under the global 
transformation of Lie groups.

\subsection{Metric derived from $SU(2)$ (I)}
Here we derive the metric from $SU(2)$ Lie group. For that purpose, 
we first give the expression of the $SU(2)$ group element in the 
form
\begin{eqnarray}
   &&U=\exp\left(\frac{i}{2} \vec{\sigma}\cdot\vec{\theta}\right)
   =\cos{\frac{|\theta|}{2}} 
   + i \frac{\vec{\sigma}  \cdot \vec{\theta}}{ |\theta| }
   \sin{\frac{|\theta|}{2}}    .
\label{e4}
\end{eqnarray}
The explicit form of the group element of $SU(3)$ is given 
in \cite{Byrd,Byrd2}.
The infinitesimal form of the group element is given by
\begin{eqnarray}
   dU=&&(-\frac{1}{2} \sin{\frac{|\theta|}{2}} 
   + i \frac{\vec{\sigma}  \cdot \vec{\theta}}{ 2|\theta| }
     \cos{\frac{|\theta|}{2}} )  d|\theta| 
\nonumber\\
   &&+i \frac{\sin(|\theta|/2)}{|\theta|} \vec{\sigma} \cdot \vec{d\theta}
   -i \frac{\vec{\sigma}  \cdot \vec{\theta}}{ |\theta|^2 }
   \sin{\frac{|\theta|}{2}} d|\theta|  .
\label{e5}
\end{eqnarray}
The expression of $U^{-1}$ is given by  
\begin{eqnarray}
   U^{-1}=\cos{\frac{|\theta|}{2}} 
   - i \frac{\vec{\sigma}  \cdot \vec{\theta}}{ |\theta| }
   \sin{\frac{|\theta|}{2}}  ,
\label{e6}
\end{eqnarray}
and we multiply this inverse factor from the left in Eq.(\ref{e5}).
Then we obtain 
\begin{eqnarray}
  &&U^{-1} dU=i \vec{\sigma} \cdot \vec{\theta}
  \left( \frac{1}{2 |\theta|}
         -\frac{\sin(|\theta|/2) \cos(|\theta|/2)}{|\theta|^{2}}\right)
    d|\theta|
\nonumber\\
 &&+i\frac{\sin(|\theta|/2) \cos(|\theta|/2)}{|\theta|}
   \sum_{a=1}^{3} \sigma_a d\theta_a
   +i\frac{\sin^{2}(|\theta|/2)}{|\theta|^2}
   \sum_{a,b,c=1}^{3} \epsilon_{abc} \theta_a d\theta_b \sigma_c  .
\label{e7}
\end{eqnarray}
After a straightforward calculation, we have 
\begin{eqnarray}
(ds)^2&&=k {\rm Tr} \left(\overline{U^{-1} dU}  U^{-1} dU \right)
\nonumber\\
&&=\frac{1}{ |\theta|^2} \left(\sum_{a} \theta_a d\theta_a \right)^2
+\frac{4 \sin^{2}(|\theta|/2) \cos^{2}(|\theta|/2)}{|\theta|^4}
\left( \sum_{a} \theta_a^2 \sum_{b} (d\theta_b)^2
 -\left( \sum_{a} \theta_a d\theta_a \right)^2 \right)
\nonumber\\
&& +\frac{4 \sin^{4}(|\theta|/2)}{|\theta|^4}
\left\{(\theta_1 d\theta_2 - \theta_2 d\theta_1)^2
      +(\theta_2 d\theta_3 - \theta_3 d\theta_2)^2
      +(\theta_3 d\theta_1 - \theta_1 d\theta_3)^2 \right\} .
\label{e8}
\end{eqnarray}
We take $k=2$ and we have the line element in the form 
\begin{eqnarray}
   (ds)^2= 2{\rm Tr} \left(\overline{U^{-1} dU}  U^{-1} dU \right)
   =\sum_{a,b} g_{a b} d\theta^a d\theta^b  .
\label{e9}
\end{eqnarray}
Then we have the expression of the metric in the form
\begin{eqnarray}
   g_{ab}= \frac{\theta_a \theta_b}{|\theta|^2}
       +\frac{4 \sin^{2}(|\theta|/2) }{|\theta|^4}
   \left( |\theta|^2 \delta_{ab} -\theta_a \theta_b \right) .
 \label{e10}
\end{eqnarray}
We put $\theta^a=\theta_a$, which is consistent with the above 
metric $g_{ab}$.
The inverse metric $ g^{ab}$, which satisfy 
$\sum_b g_{ab} g^{bc}=\delta^c_{a}$, is given by 
\begin{eqnarray}
   g^{ab}= \frac{\theta^a \theta^b}{|\theta|^2}
       +\frac{1}{4 \sin^{2}(|\theta|/2) }
   \left( |\theta|^2 \delta^{ab} -\theta^a \theta^b \right) .
\label{e11}
\end{eqnarray}

This metric $g_{ab}$ satisfy the Einstein equation 
with the cosmological constant in the form
\begin{eqnarray}
   R_{ab}-\frac{1}{2} R g_{ab}+ \Lambda g_{ab}=0 ,
\label{e12}
\end{eqnarray}
with $\Lambda=1/4$ in this case.
This is the Einstein metric, because it satisfies the 
equation
\begin{eqnarray}
   R_{ab}=2 \Lambda g_{ab} .
\label{e13}
\end{eqnarray} 

\subsection{Metric derived from $SU(2)$(II):\quad
Euler angle parametrization}
We parametrize the $SU(2)$ group element with the Euler angle 
in the form
\begin{eqnarray}
   && U=U_z(\phi) U_x(\theta) U_z(\psi) \nonumber\\
   &&=
   \left( \begin{array}{cc} e^{i \phi/2} & 0\\
                              0 & e^{- i \phi/2}
            \end{array}
   \right)
   \left( \begin{array}{cc} \cos(\theta/2) & i \sin(\theta/2)\\
                            i \sin(\theta/2)  & \cos(\theta/2)
          \end{array}   
   \right)
   \left( \begin{array}{cc} e^{i \psi/2} & 0\\
                              0 & e^{- i \psi/2}
            \end{array}
   \right)
\nonumber\\
  &&=\left( \begin{array}{cc} e^{i (\phi+\psi)/2}\cos(\theta/2) & 
    i e^{i (\phi-\psi)/2}\sin(\theta/2)\\
    i e^{-i (\phi-\psi)/2}\sin(\theta/2) & 
    e^{-i (\phi+\psi)/2}\cos(\theta/2)
            \end{array}
   \right)    .
\label{e14}
\end{eqnarray} 
From this expression, we have 
\begin{eqnarray}
   U^{-1}dU=
   \left( \begin{array}{cc} 
      \frac{i}{2} \cos(\theta) d\phi + \frac{i}{2} d\psi
     & \frac{i}{2} e^{-i \psi}d\theta-\frac{1}{2} e^{-i \psi} 
       \sin(\theta) d\phi \\
       \frac{i}{2} e^{i \psi}d\theta+\frac{1}{2} e^{i \psi} 
       \sin(\theta) d\phi               
     & -\frac{i}{2} \cos(\theta) d\phi - \frac{i}{2} d\psi
            \end{array}
   \right)  .
\label{e15}
\end{eqnarray} 
Then we define the line element in the form 
\begin{eqnarray}
   (ds)^2&&=2 {\rm Tr} \left(\overline{U^{-1} dU}  U^{-1} dU \right)
\nonumber\\
   &&=(\cos(\theta) d\phi + d\psi)^2+(d\theta)^2+\sin^2(\theta) (d\phi)^2
\label{e16}\\
   &&=(d\theta)^2+(d\psi)^2+(d\phi)^2+2 \cos(\theta) d\psi d\phi .
\label{e17}
\end{eqnarray} 
This is the standard expression of the line element 
for $SU(2)$ \cite{Gubser}.
If we take $\theta^1=\theta$, $\theta^2=\phi$, $\theta^3=\psi$,
non-zero elements of the metric are given by
\begin{eqnarray}
 && g_{11}=1, \quad g_{22}=1, \quad g_{33}=1, 
    \quad g_{23}=g_{32}=\cos(\theta^1) , 
\label{e18}\\
 && g^{11}=1, \quad g^{22}=\frac{1}{\sin^2(\theta^1)}, 
    \quad g^{33}=\frac{1}{\sin^2(\theta^1)}, 
    \quad g^{23}=g^{32}=-\frac{\cos(\theta^1)}{\sin^2(\theta^1)} .
\label{e19}
\end{eqnarray}
In this metric, we can easily see that this line element 
$(ds)^2$ is invariant under 
\begin{eqnarray}
 &&i) \quad \phi \rightarrow \phi^{'}=\phi+\xi, \quad 
     \theta \rightarrow \theta^{'}=\theta, \quad 
     \psi \rightarrow \psi^{'}=\psi , 
\label{e20} \\
&&ii) \quad \psi \rightarrow \psi^{'}=\psi+\xi, \quad 
     \theta \rightarrow \theta^{'}=\theta, \quad 
     \phi \rightarrow \phi^{'}=\phi , 
\label{e21} 
\end{eqnarray}
with global parameter $\xi$. The third transformation, which 
makes $(ds)^2$ invariant, is quite complicated. 

The above expression of the metric $g_{ab}$ also becomes 
the Einstein metric and it satisfies Eq.(\ref{e13}). 

\section{ Metric derived from the quadratic invariant manifold for the 
classical Lie groups}
Classical Lie groups such as $SU(N)$, $SO(N)$, $Sp(N)$,
are defined in such a way as the set of transformations 
to make the special quadratic form invariant.
For example, $SO(N)$ Lie group is defined in such a way 
as 
\begin{eqnarray}
\sum_{i=1}^{N} x_i^{2}=({\rm invariant}) ,
\label{e22}
\end{eqnarray}
under the transformation 
\begin{eqnarray}
\left( \begin{array}{c} x^{'}_1\\ x^{'}_2 \\ \cdots 
       \\x^{'}_{N}
       \end{array}
   \right)
=U \left( \begin{array}{c} x_1\\ x_2 \\ \cdots \\x_{N}
       \end{array}
   \right) .
\label{e23}
\end{eqnarray}
Then the metric derived on this quadratic invariant manifold 
have the relation with the metric derived from the classical 
Lie groups. The symmetry of the metric derived from the quadratic 
invariant manifold is the subset of the symmetry of the metric 
derived from the associated classical Lie groups.
For $SO(N)$ case, we parametrize the above quadratic 
invariant form with $N-1$ angle variables 
in the form $ x_i=x_i(\theta^1, \theta^2, \cdots, \theta^{N-1})$.
We define the vector $B_{i, a}$ in the form 
\begin{eqnarray}
  B_{i, a}=\frac{\partial x_i}{\partial \theta^a} . 
\label{e24}
\end{eqnarray}
Then the metric on the quadratic invariant manifold is given by
\begin{eqnarray}
   g_{ab}=\sum_{i=1}^{N} B_{i, a} B_{i, b} .
\label{e25}
\end{eqnarray}

For the $SO(3)$ case, we parametrize the invariant quadratic manifold 
\begin{eqnarray}
   x_1^{2}+x_2^{2}+x_3^{2}=1  ,
\label{e26}
\end{eqnarray}
in the form 
\begin{eqnarray}
   &&  x_1=\sin{\theta}\cos{\phi} ,
   \quad  x_2=\sin{\theta}\sin{\phi} , \quad
   x_3=\cos{\theta} .
   \label{e27}
\end{eqnarray}
Then the vector $\vec{B}_a$ given by 
\begin{eqnarray}
  &&\vec{B}_1 =\frac{\partial x_i}{\partial \theta}=
  (\cos{\theta} \cos{\phi}, \cos{\theta} \sin{\phi}, 
  -\sin{\theta}) , 
\label{e28}\\
  &&\vec{B}_2 =\frac{\partial x_i}{\partial \phi} =
  (-\sin{\theta} \sin{\phi}, \sin{\theta} \cos{\phi}, 0) .
\label{e29}
\end{eqnarray}
Then the metric is given by 
\begin{eqnarray}
  g_{11}=<\vec{B}_1,\vec{B}_1>=1, 
  \quad g_{22}=<\vec{B}_2,\vec{B}_2>=\sin^2{\theta}, \quad 
  g_{12}=g_{21}=<\vec{B}_1,\vec{B}_2>=0   .
\label{e30}
\end{eqnarray}
The invariant length is given by
\begin{eqnarray}
  (ds)^2=(d \theta)^2+\sin^2{\theta}(d \phi)^2 .
\label{e31}
\end{eqnarray}
If we compare Eq.(\ref{e31}) with Eq.(\ref{e16}), we can see that 
the symmetry of the line element Eq.(\ref{e31}) is the subset of the 
symmetry of the line element of Eq.(\ref{e16}).

\section{Summary}
We give general form of the expression of the metric derived 
from Lie groups. Then we give the explicit expression of the metric 
derived from $SU(2)$ Lie group in two different ways. We compare 
the metric on the manifold $S^2$ with the metric derived from 
$SU(2)$ Lie group. From that observation, we conjecture that the 
metric becomes the Einstein metric 
$R_{\mu \nu}=2\Lambda g_{\mu \nu}$ if the metric is derived from 
classical Lie groups $SO(n)$, $SU(n)$, $Sp(n)$.

\vspace{20mm}

\end{document}